\documentclass[reqno]{amsart}
\usepackage{amsmath,amssymb,stmaryrd}
\newtheorem{thm}{Theorem}[section]
\newtheorem{lem}{Lemma}[section]

\newtheorem{propdefi}{Proposition-Definition}[section]

\def \rn {\mathbb{R}^4}
\def \rr {\mathbb{R}}
\def \nn {\mathbb{N}}
\def \tu {\tilde{u}}
\def \wk {w_k}
\def \xk {x_k}
\def \rk {r_k}

\def \yk {y_k}

\def \Vk {V_k}
\def \tVk {\tilde{V}_k}
\def \uk {u_k}

\def \mk {\mu_k}
\def \vk {v_k}
\def \tuk {\tilde{u}_k}
\def \vphi {\varphi}

\def \beq {\begin{eqnarray*}}
\def \eeq {\end{eqnarray*}}
\def \beqn {\begin{eqnarray}}
\def \eeqn {\end{eqnarray}}
\def \bequa {\begin{equation}}
\def \eequa {\end{equation}}

\title[Fourth order equation]{Concentration phenomena for a fourth order equations with exponential growth:\\
the radial case}
\author{Fr\'ed\'eric Robert}
\email{frobert@math.unice.fr}
\address{Universit\'e de Nice-Sophia Antipolis, Laboratoire J.A.Dieudonn\'e, Parc Valrose, 06108 NICE Cedex 2, France.}
\date{December, 7th 2005.}
\begin{document}
\begin{abstract}
We let $\Omega$ be a smooth bounded domain of $\rn$ and a sequence of fonctions $(\Vk)_{k\in\nn}\in C^0(\Omega)$ such that $\lim_{k\to +\infty}\Vk=1$ in $C^0_{loc}(\Omega)$. We consider a sequence of functions $(\uk)_{k\in\nn}\in C^4(\Omega)$ such that $$\Delta^2\uk=\Vk e^{4\uk}$$
in $\Omega$ for all $k\in\nn$. We address in this paper the question of the asymptotic behaviour of the $(\uk)'s$ when $k\to +\infty$. The corresponding problem in dimension 2 was considered by Br\'ezis-Merle and Li-Shafrir (among others), where a blow-up phenomenon was described and where a quantization of this blow-up was proved. Surprisingly, as shown by Adimurthi, Struwe and the author in \cite{ars}, a similar quantization phenomenon does not hold for this fourth order problem. Assuming that the $\uk$'s are radially symmetrical, we push further the analysis of \cite{ars}. We prove that there are exactly three types of blow-up and we describe each type in a very detailed way.
\end{abstract}
\maketitle

\section{Introduction}
Let $\Omega$ be a bounded domain of $\rn$. Let $(\Vk)_{k\in\nn}\in C^0(\Omega)$ be a sequence such that
\begin{equation}\label{lim:Vk}
\lim_{k\to +\infty}\Vk=1
\end{equation}
in $C^0_{loc}(\Omega)$. Let $(\uk)_{k\in\nn}$ be a sequence of functions in $C^4(\Omega)$ such that
$$\Delta^2\uk= \Vk e^{4\uk}\eqno{(E)}$$
in $\Omega$ for all $k\in\nn$. Here and in the sequel, $\Delta=-\sum\partial_{ii}$ is the Laplacian with minus sign convention. 
In this paper, we address the question of the asymptotics of the $\uk$'s when $k\to +\infty$. A natural (and simple) behavior is when 
there exists
$u\in C^3(\Omega)$ such that, up to a subsequence, 
\begin{equation}\label{eq:cv}
\lim_{k\to +\infty}\uk=u
\end{equation}
in $C^{3}_{loc}(\Omega)$. In this situation, we say that $(\uk)_{k\in\nn}$ is relatively compact in $C^3_{loc}(\Omega)$. 
However, the structure of equation $(E)$ is much richer due to its scaling invariance properties. 
The scaling invariance is as follows. Given $k\in\nn$, $x_k\in\Omega$
and $\mu_k>0$ , we let
\begin{equation}\label{eq:inv:uk}
\tilde{u}_{k}(x):=\uk(\xk+\mu_k x)+\ln\mu_k
\end{equation}
for all $x\in\mu_k^{-1}(\Omega-\xk)$. Letting $\tilde{V}_k(x)=\Vk(x_k+\mu_k x)$ for all $x\in \mu_k^{-1}(\Omega-\xk)$, we get 
that the rescaled function $\tilde{u}_{k}$ satisfies
$$\Delta^2\tilde{u}_{k}=\tilde{V}_k e^{4\tilde{u}_{k}}$$
on $\mu_k^{-1}(\Omega-\xk)$ -- an equation like $(E)$. This scaling invariance forces some situations more subtle than (\ref{eq:cv}) to
happen. A very basic example is the following: we consider a sequence $(\mu_k)_{k\in\nn}\in\rr_{>0}$ such that $\lim_{k\to +\infty}\mu_k=0$. Let a function $v\in C^4(\rn)$ such that $e^{4v}\in L^1(\rn)$ and
\bequa\label{eq:lim}
\Delta^2 v=e^{4v}.
\eequa
The simplest example is the function $x\mapsto\ln\frac{\sqrt{96}}{\sqrt{96}+|x|^2}$. For any $k\in\nn$, we define the function
$$f_k(x)=v\left(\mk^{-1}x\right)-\ln\mk$$
for all $x\in\rn$. Then $f_k$ satisfies $(E)$ with $\Vk\equiv 1$ for all $k\in\nn$, but the sequence 
$(f_k)_{k\in\nn}$ does not converge in $C^{0}_{loc}(\rn)$: indeed, we have that
$$\lim_{k\to +\infty}f_k(0)=+\infty\hbox{ and }\Vk e^{4 f_k}\, dx\rightharpoonup \left(\int_{\rn}e^{4v}\, dx\right)\delta_0$$
when $k\to +\infty$ weakly for the convergence of measures. Here and in the sequel, $\delta_0$ denotes the Dirac mass at $0$, and we say that the energy of the sequence $(f_k)$ is $\int_{\rn}e^{4v}\, dx$. Scaling as in (\ref{eq:inv:uk}), we get that
$$\lim_{k\to +\infty}f_k(\mu_k x)+\ln\mk=v(x)$$
for all $x\in\rn$. In other words, $(f_k)$ converges to $v$ up to rescaling. Concerning terminology, we say that the sequence $(\uk)_{k\in\nn}$ blows-up if it is not relatively 
compact in $C^3_{loc}(\Omega)$, so that, up to any subsequence, (\ref{eq:cv}) does not hold. In the above 
example, the $(f_k)$'s blow up. In this paper,
we are concerned with the blow-up behavior of solutions of $(E)$.

\medskip\noindent In dimension two, the corresponding problem involves the Laplacian (and not the bi-Laplacian). This problem has been
studied (among others) by Br\'ezis-Merle \cite{brezismerle} and
Li-Shafrir \cite{lishafrir}. We also refer to Druet \cite{druet} and
Adimurthi-Struwe \cite{adistruwe} for the description of equations with more intricate nonlinearities and to Tarantello \cite{t} for equations with singularities. An important phenomenon that holds in dimension two is the quantization of the energy. Following standard terminology, we say that there is quantization if there exists a positive constant $C_m>0$ such that the energy of any blowing-up sequence of solutions to the equation under consideration is (roughly speaking) asymptotically a multiple of $C_m$. In particular, when blow-up occurs, the sequence of solutions carries at least the energy $C_m$ or carries no energy.

\medskip\noindent Surprisingly, such a quantization result is false when we come back to our initial four-dimensional problem $(E)$. Let $\lambda\in(0,+\infty)$ arbitrary: in a joint work with Adimurthi and Michael Struwe \cite{ars}, we exhibit a sequence of solutions to $(E)$ that blows-up, carries the energy $\lambda$ and develop singularities on a $3-$dimensional hypersurface of $\rn$. Still in \cite{ars}, we described the behaviour of arbitrary solutions to $(E)$ and proved that any blowing-up sequence $(\uk)_{k\in\nn}$ concentrates at the zero set of a nonpositive nontrivial bi-harmonic function, and that outside this set, $\lim_{k\to +\infty}\uk=-\infty$ uniformly. In view of the results of \cite{ars}, giving a more precise description requires additional hypothesis on $(\uk)_{k\in\nn}$. 

\medskip\noindent A natural hypothesis is to impose a Navier boundary condition, (that is $\uk=\Delta\uk=0$ on $\partial\Omega$) or a Dirichlet boundary condition (that is $\uk=\frac{\partial\uk}{\partial\nu}=0$ on $\partial\Omega$): actually, in these cases, we get that there is no blow-up and we recover relative compactness, these claims are easy consequences of the result in \cite{wei}. Wei \cite{wei} also studied the case where $\Delta\uk=0$ on $\partial\Omega$ and $\uk=c_k$ on $\partial\Omega$, where $(c_k)_{k\in\nn}\in\rr$ is a sequence of real numbers such that $\lim_{k\to +\infty}c_k=-\infty$: in this context, Wei described precisely the asymptotics and recovered a quantization result as in Li-Shafrir. In \cite{r1}, we consider the case where the $L^1-$norm of $\Delta\uk$ is uniformly bounded on a given subset of $\Omega$: in this context, we also recover a quantization result (that is the energy of a blowing-up solution is a multiple of an explicit constant).

\medskip\noindent In the present paper, we consider the case when $\Omega=B$ is a ball and when the $\uk$'s are radially symmetrical with respect to the center of the ball for all $k\in\nn$. Without loss of generality, we assume that $B=B_1(0)$ is the unit ball of $\rn$ centered at $0$. In this rather natural situation, and contrary to the situation considered in \cite{r1}, there is no quantization. This phenomenon is due to the abundance of solutions to equation \eqref{eq:lim} (see C.S.Lin \cite{lin}), contrary to the two-dimensional corresponding equation, where up to affine transformations, there is only one solution. 



\medskip\noindent Let $(\uk)_{k\in\nn}$ be a sequence of blowing-up solutions to $(E)$, with a sequence $(\Vk)_{k\in\nn}\in C^0(B)$ such that \eqref{lim:Vk} holds. Assuming that the $\uk$'s are radially symmetrical, the first step in studying the blow-up behavior of the $(\uk)$'s is to prove that $\Vk e^{4\uk}\, dx$ converges to the product of a real number (refered to as the energy) by a Dirac mass at $0$ for the convergence of measures when $k\to +\infty$: it is much more tricky to have informations about the energy in front of the Dirac mass, and this is the object of Theorem \ref{thm:main}. The intricate issue in this theorem concerns the localization of the energy at the microscopic level. More precisely, after rescaling as in \eqref{eq:inv:uk}, we prove (in general) that the $(\uk)$'s converge when $k\to +\infty$ to a solution $v\in C^4(\rn)$ of \eqref{eq:lim} such that $e^{4v}\in L^1(\rn)$: since the $L^1-$norm is invariant under the rescaling \eqref{eq:inv:uk}, we get that there exists a sequence $(\rk)_{k\in\nn}$ of positive real numbers such that $\lim_{k\to +\infty}\rk=0$ and such that the $L^1-$norm of $e^{4\uk}$ in $B_{\rk}(0)$ converges to the $L^1-$norm of $e^{4v}$ in $\rn$. The difficult step is to prove that there is no energy left outside this ball of radius $\rk$ when $k\to +\infty$, and so, in other words, that the $L^1-$norm of $e^{4\uk}$ outside $B_{\rk}(0)$ goes to $0$ when $k\to +\infty$. Refering to standard terminology, this corresponds to provinge that there is no energy lost in the necks. Our main result is the following. 


\begin{thm}\label{thm:main} Let $(\Vk)_{k\in\nn}\in C^0(B)$ be a sequence of functions such that (\ref{lim:Vk}) holds. Let $(\uk)_{k\in\nn}$ be a family of functions in $C^4(B)$ which are solutions to $(E)$. We assume that there exists $\Lambda\in\rr$ such that 
$$\int_B \Vk e^{4\uk}\, dx\leq \Lambda$$
for all $k\in\nn$ and that the $(\uk)$'s blow-up, that is \eqref{eq:cv} does not hold for any subsequence. In addition, we assume that $\uk$ is radially symmetrical for all $k\in\nn$. Then, up to a subsequence, there exists $\alpha\in [0,16\pi^2]$ such that
$$\Vk e^{4\uk}\, dx\rightharpoonup \alpha\delta_0$$
when $k\to +\infty$ for the convergence of measures. More precisely,\par

\smallskip (i) either there exists $C>0$ such that, up to a subsequence, $\uk(0)\leq C$ for all $k\in\nn$: then $\alpha=0$ and $\lim_{k\to +\infty}\uk=-\infty$ uniformly locally on $B\setminus\{0\}$

\smallskip (ii) or $\lim_{k\to +\infty}\uk(0)=+\infty$. In this situation, for any $\delta\in (0,1)$, we have that
$$\lim_{R\to +\infty}\lim_{k\to +\infty}\int_{B_\delta(0)\setminus B_{R e^{-\uk(0)}}(0)}\Vk e^{4\uk}\, dx=0.$$
In addition, still in case (ii), the asymptotic behavior at the scale $e^{-\uk(0)}$ is ruled as follows:

\smallskip\noindent(ii.a) if $\alpha=16\pi^2$, then
$$\lim_{k\to +\infty}\left(\uk(e^{-\uk(0)}x)-\uk(0)\right)= \ln\frac{\sqrt{96}}{\sqrt{96}+|x|^2}$$
for all $x\in\rn$. Moreover, this convergence holds in $C_{loc}^3(\rn)$.

\smallskip\noindent(ii.b) if $\alpha\in (0,16\pi^2)$, then there exists $v\in C^4(\rn)$ such that $e^{4v}\in L^1(\rn)$ and
$$\lim_{k\to +\infty}\left(\uk(e^{-\uk(0)}x)-\uk(0)\right)= v(x)$$
when $k\to +\infty$ for all $x\in\rn$. Moreover, this convergence holds in $C_{loc}^3(\rn)$ and there exists $\lambda>0$ such that $\lim_{|x|\to +\infty}\frac{v(x)}{|x|^2}=-\lambda$.

\smallskip\noindent(ii.c) If $\alpha=0$, then $\lim_{k\to +\infty}e^{-2\uk(0)}\Delta\uk(0)=+\infty$ and we have that
$$\lim_{k\to+\infty}\frac{\uk(e^{-\uk(0)}x)-\uk(0)}{e^{-2\uk(0)}\Delta\uk(0)}\to -\frac{|x|^2}{8}$$
when $k\to +\infty$ for all $x\in\rn$. Moreover, this convergence holds in $C_{loc}^3(\rn)$.
\end{thm}
\noindent Note that this theorem is optimal: for any $\alpha\in [0,16\pi^2]$, we exhibit in section \ref{sec:examples} examples of blowing-up solutions to $(E)$ such that their energy converges to $\alpha$. Note also that this theorem is specific to the radial case and does not hold in general for nonradial solutions (see for instance Adimurthi-Robert-Struwe \cite{ars}).

\medskip\noindent In a joint work with Olivier Druet \cite{druetrobert}, we
studied the corresponding problem on four-dimensional Riemannian
manifolds, where the bi-Laplacian is replaced by a fourth-order
elliptic operator refered to as $P$: when the kernel of $P$ is such
that $\hbox{Ker } P=\{constants\}$, we get that blow-up occurs at finitely many isolated points, and that each point carries exactly the energy $16\pi^2$. Note that in the context of Theorem \ref{thm:main}, the kernel of the bi-Laplacian contains more than the constant functions. Related references in the context of Riemannian
manifolds are Malchiodi \cite{malchiodi} and Malchiodi-Struwe \cite{malchiodistruwe}. As a
remark, the corresponding question in dimension $n\geq 5$ was
considered in Hebey-Robert \cite{hebeyrobert}, we refer also to Hebey-Robert-Wen \cite{hrw}.

\medskip\noindent This paper is organized as follows. In section \ref{sec:examples}, we exhibit examples of blowing-up solutions to $(E)$ having any energy ranging in $[0,16\pi^2]$. In sections \ref{sec:prem} to \ref{sec:lim0}, we prove Theorem \ref{thm:main}. More precisely, in section \ref{sec:prem}, we introduce the three types of convergence that correspond to the cases $\alpha=16\pi^2$, $0<\alpha<16\pi^2$ and $\alpha=0$ of Theorem \ref{thm:main}. The case $\alpha=0$ of Theorem \ref{thm:main} is proved in section \ref{sec:liminfty}. The case $0<\alpha<16\pi^2$ of Theorem \ref{thm:main} is proved in section \ref{sec:limK}. The case $\alpha=16\pi^2$ of Theorem \ref{thm:main} is proved in section \ref{sec:lim0}. In the sequel, $C$ denotes a positive constant, with value allowed to change from one line to the other. Note also that all the convergence results are up to a subsequence, even when it is not precised.

\medskip\noindent{\bf Acknowledgements:} the author thanks Adimurthi and Michael Struwe for having suggested him to work on these questions and for stimulating discussions. The author is indebted to Olivier Druet and Emmanuel Hebey for stimulating discussions.

\section{Examples of Log- and Quadratic-Convergences}\label{sec:examples}

We exhibit situations in which the three patterns {\it (ii.a)}, {\it .b)} and {\it .c)} occur. 

\subsection{Log-Convergence}
We let $k\in\nn^\star$ and define the function
$$\uk(x):=\ln \frac{k\sqrt{96}}{\sqrt{96}+k^2|x|^2}$$
for all $x\in\rn$. We have that 
$$\uk(e^{-\uk(0)} x)-\uk(0)=\ln \frac{\sqrt{96}}{\sqrt{96}+|x|^2}$$
for all $x\in\rn$. As easily checked,
$$\Delta^2\uk=e^{4\uk}$$
and 
$$\Vk e^{4\uk}\rightharpoonup 16\pi^2\delta_0$$
in $B_1(0)$ in the sense of measures when $k\to +\infty$, and we are in the situation described in {\it (ii.a)}.

\subsection{Quadratic-Convergence (I)} 
We let $\alpha\in (0,16\pi^2)$. It follows from \cite{changchen} that there exists $v\in C^4(\rn)$ such that $v\leq v(0)=0$ and $\Delta^2 v= e^{4v}$ in $\rn$ and $\int_{\rn}  e^{4v}\, dx=\alpha$. For any $k\in\nn^\star$, we define the function
$$\uk(x):=v\left(k x\right)+\ln k$$
for all $x\in\rn$. We have that 
$$\Delta^2\uk=e^{4\uk}$$
in $B$ for all $k\in\nn^\star$ and we have that
$$\uk(e^{-\uk(0)}x)-\uk(0)=v(x)$$
for all $x\in\rn$ and all $k\in\nn^\star$. 
In addition, we have that
$$\Vk e^{4\uk}\, dx\rightharpoonup \alpha\delta_0$$
in $B_1(0)$ in the sense of measures when $k\to +\infty$ and we are in the situation described in {\it (ii.b)}.

\subsection{Quadratic-Convergence (II)}
We let the unique radially symmetrical function $\vphi\in C^4(\rn)$ such that $\Delta^2\vphi=e^{-\frac{|x|^2}{2}}$ in $\rn$, $\vphi(0)=\Delta\vphi(0)=0$.
We let 
$$\uk(x):=\ln k-\frac{k^6|x|^2}{8}+k^{-8}\vphi\left(k^3 x\right)$$
for all $k\in\nn^\star$ and $x\in\rn$. We define
$$\Vk=e^{-4\uk}\Delta^2\uk$$
for all $k\in\nn^\star$. All these functions are explicit (see \cite{ars}) and we get that
$$\lim_{k\to +\infty}\Vk=1\hbox{ in }C^0_{loc}(\rn).$$
Moreover, we have that
$$\Vk e^{4\uk}\, dx\rightharpoonup 0$$
in $B_1(0)$ in the sense of measures when $k\to +\infty$ and we are in the situation described in {\it (ii.c)}. We refer to \cite{ars} for details about these assertions. A similar method permits to construct families $(\uk)$ and $(\Vk)$ such that \ref{lim:Vk} and $(E)$ hold and such that $\uk\leq\uk(0)=0$ and $\Vk e^{4\uk}\, dx\rightharpoonup 0$ when $k\to +\infty$, and we are in the situation described in {\it (i)}.


\medskip\noindent These three examples show that for any $\alpha\in [0,16\pi^2]$, their exists a blowing-up sequence of solutions to $(E)$ with energy $\alpha$.

\section{Preliminary estimates for $(E)$}\label{sec:prem}
We let $B$ be the open unit ball of $\rn$ and $(\Vk)_{k\in\nn}\in C^0(B)$ a sequence such that
\bequa\label{hyp:Vk}
\lim_{k\to +\infty}\Vk=1\hbox{ in }C^0_{loc}(B).
\eequa
We let $(\uk)_{k\in\nn}\in C^4(B)$ such that for any $k\in\nn$, we have that
\bequa\label{eq:uk}
\Delta^2\uk=\Vk e^{4\uk}
\eequa
in $B$. We assume that is there exists $\Lambda\in\rr$ such that for any $k\in\nn$, we have that
\bequa\label{bnd:nrj}
\int_B  e^{4\uk}\, dx\leq \Lambda.
\eequa
We assume that the function $\uk$ is radially symmetrical with respect to the center of the unit ball $B$, that is $0$. For any radially symmetrical function $h$, there exists $\tilde{h}$ defined on an interval of $[0,+\infty)$ such that $h(x)=\tilde{h}(|x|)$ for all $x$ such that this expression makes sense. With a standard abuse of notation, we write $h(r)$, $h'(r)$, etc for $\tilde{h}(r)$, $\tilde{h}'(r)$ respectively. This section is devoted to the proof of general estimates on the $(\uk)$'s and to the definition of the three types of convergnce that will let us distinguish the three situations of blow-up in Theorem \ref{thm:main}.

\medskip\noindent{\bf Step \ref{sec:prem}.1:} We first deal with the behaviour of $\uk$ on subsets where it is bounded from above:

\begin{lem}\label{lem:wpe} Let $(\Vk)_{k\in\nn}\in C^0(B)$ and $(\uk)_{k\in\nn}\in C^4(B)$ such that (\ref{hyp:Vk}), (\ref{eq:uk}) and (\ref{bnd:nrj}) hold. We assume that $\uk$ is radially symmetrical for all $k\in\nn$. We let $\omega\subset\subset B$. Then, there exists $C(\omega)>0$ such that
\bequa\label{estim:wpe}
|x|e^{\uk(x)}\leq C(\omega)
\eequa
for all $x\in\omega$ and all $k\in\nn$.
\end{lem}

\smallskip\noindent{\it Proof of Lemma \ref{lem:wpe}:} We let $\delta_1\in (0,1)$ such that $\omega\subset B_{\delta_1}(0)$. We let $\delta_2\in (\delta_1,1)$. Since (\ref{hyp:Vk}) holds, we assume without loss of generality that 
\bequa\label{inf:Vk}
\Vk(x)\geq\frac{1}{2}
\eequa
for all $x\in B_{\delta_2}(0)$ and all $k\in\nn$. With equation (\ref{eq:uk}), we get that $\Delta(\Delta\uk)>0$ on $B$, and then $\Delta\uk$ (considered as a function of $r\in [0,1)$) is strictly decreasing on $[0,\delta_2]$. We distinguish three situations:

\smallskip\noindent{\it Case \ref{sec:prem}.1.1:} We assume that $\Delta\uk\geq 0$ on $B_{\delta_2}(0)$. In this situation, we get that $\uk$ is decreasing on $[0,\delta_2]$. We let $x\in B_{\delta_1}(0)$. With (\ref{bnd:nrj}) and (\ref{inf:Vk}), we get that
$$\Lambda\geq \int_{B_{|x|}(0)}\Vk e^{4\uk}\, dy\geq \frac{e^{4\uk(x)}}{2}\hbox{Vol}(B_{|x|}(0))\geq \frac{\pi^2|x|^4 e^{4\uk(x)}}{4}.$$
In particular, \eqref{estim:wpe} holds in Case \ref{sec:prem}.1.1.\par


\smallskip\noindent{\it Case \ref{sec:prem}.1.2:} We assume that $\Delta\uk\leq 0$ on $B_{\delta_2}(0)$. In this situation, we get that $\uk$ is increasing on $[0,\delta_2]$. We let $x\in B_{\delta_1}(0)$. With (\ref{bnd:nrj}) and (\ref{inf:Vk}), we get that
\beq
\Lambda & \geq & \int_{B_{\delta_2}(0)\setminus B_{|x|}(0)}\Vk e^{4\uk}\, dy\geq \frac{e^{4\uk(x)}}{2}\hbox{Vol}(B_{\delta_2}(0)\setminus B_{|x|}(0))\\
&\geq &\frac{\pi^2(\delta_2^4-|x|^4) e^{4\uk(x)}}{4}\geq \frac{\pi^2(\delta_2^4-\delta_1^4) e^{4\uk(x)}}{4}.
\eeq
In particular, \eqref{estim:wpe} holds in Case \ref{sec:prem}.1.2.\par


\smallskip\noindent{\it Case \ref{sec:prem}.1.3:} We assume that $\Delta\uk$ takes some positive and some negative values in $B_{\delta_2}(0)$. Since $\Delta\uk$ is decreasing, there exists $s_k\in (0,\delta_2)$ such that 
$$\Delta\uk>0\hbox{ in }[0,s_k),\;  \Delta\uk(s_k)=0,\hbox{ and }\Delta\uk<0\hbox{ in }(s_k,\delta_2].$$
In particular, there exists $\tau_k\in [s_k,\delta_2]$ such that $\uk$ is decreasing in $[0,\tau_k)$ and $\uk$ is increasing in $[\tau_k,\delta_2]$ (note that the case $\tau_k=\delta_2$ is possible). We let $x\in B_{\delta_1}(0)$. If $|x|\leq \tau_k$, we proceed as in Case \ref{sec:prem}.1.1. If $|x|\geq \tau_k$, we proceed as in \ref{sec:prem}.1.2. In particular, \eqref{estim:wpe} holds in Case \ref{sec:prem}.1.2.\par


\noindent These three cases prove Lemma \ref{lem:wpe}. \hfill$\Box$

\medskip\noindent{\bf Step \ref{sec:prem}.2:} The preceding step permits us to deal with the convergence outside $0$. This is the object of the following Lemma:

\begin{lem}\label{cv:out0}
Let $(\Vk)_{k\in\nn}\in C^0(B)$ and $(\uk)_{k\in\nn}\in C^4(B)$ such that (\ref{hyp:Vk}), (\ref{eq:uk}) and (\ref{bnd:nrj}) hold. We assume that $\uk$ is radially symmetrical for all $k\in\nn$. Then we are in one and only one of the following situations:

(a) there exists $u\in C^4(B\setminus\{0\})$ such that, up to a subsequence,
$$\lim_{k\to +\infty}\uk=u\hbox{ in }C^3_{loc}(B\setminus\{0\}).$$

(b) there exists a sequence $(a_k)_{k\in\nn}\in\rr_{>0}$ such that $\lim_{k\to +\infty}a_k=+\infty$, there exists $\varphi\in C^4(B\setminus\{0\})$ such that $\Delta^2\varphi=0$, $\varphi<0$, and such that
$$\lim_{k\to +\infty}\frac{\uk}{a_k}=\varphi\hbox{ in }C^3_{loc}(B\setminus\{0\}).$$ 
In particular, $\uk\to-\infty$ uniformly on every compact subset of $B\setminus\{0\}$.
\end{lem}

\noindent We omit the proof of the Lemma: it is a direct consequence of the results of \cite{ars} combined with Lemma \ref{lem:wpe}. We refer to \cite{ars} for details.

\medskip\noindent{\bf Step \ref{sec:prem}.3:} This short step is devoted to the case when $\uk$ is bounded from above. More precisely we have:

\begin{lem}\label{cv:uk:bnd}
Let $(\Vk)_{k\in\nn}\in C^0(B)$ and $(\uk)_{k\in\nn}\in C^4(B)$ such that (\ref{hyp:Vk}), (\ref{eq:uk}) and (\ref{bnd:nrj}) hold. We assume that $\uk$ is radially symmetrical for all $k\in\nn$. We assume that there exists $\delta_0\in (0,1)$ and $C(\delta_0)>0$ such that 
\bequa\label{hyp:bnd:sup}
\uk(x)\leq C(\delta_0)
\eequa
for all $x\in B_{\delta_0}(0)$. Then we are in one and only one of the following situations:

\smallskip(a) there exists $u\in C^4(B)$ such that, up to a subsequence,
$$\lim_{k\to +\infty}\uk=u\hbox{ in }C^3_{loc}(B).$$
In particular, $$\Vk e^{4\uk}\, dx\rightharpoonup e^{4u}\, dx$$
when $k\to +\infty$ in the sense of measures.

\smallskip(b) there exists a sequence $(a_k)_{k\in\nn}\in\rr_{>0}$ such that $\lim_{k\to +\infty}a_k=+\infty$, there exists $\varphi\in C^4(B)$ such that $\Delta^2\varphi=0$, $\varphi<0$ in $B\setminus\{0\}$, and such that
$$\lim_{k\to +\infty}\frac{\uk}{a_k}=\varphi\hbox{ in }C^3_{loc}(B).$$ 
In particular, $$\Vk e^{4\uk}\, dx\rightharpoonup 0$$
when $k\to +\infty$ in the sense of measures.
\end{lem}

\medskip\noindent{\it Proof of Lemma \ref{cv:uk:bnd}:} it follows from (\ref{estim:wpe}) and (\ref{hyp:bnd:sup}) that for any $\delta\in (0,1)$, there exists $C(\delta)>0$ such that $\uk(x)\leq C(\delta)$ for all $x\in\overline{B}_\delta(0)$. We proceed as in \cite{ars} and we obtain that the function $\varphi$ in Lemma \ref{cv:out0} is defined on the whole domain $B$, and is bi-harmonic in $B$. Since $\varphi$ is radially symmetrical, $\varphi\leq 0$ and $\varphi\not\equiv 0$, we get that $\varphi<0$ in $B\setminus\{0\}$. This proves Lemma \ref{cv:uk:bnd}.
\hfill$\Box$

\medskip\noindent{\bf Step \ref{sec:prem}.4:} 
\begin{lem}\label{lem:uk:infty}
Let $(\Vk)_{k\in\nn}\in C^0(B)$ and $(\uk)_{k\in\nn}\in C^4(B)$ such that (\ref{hyp:Vk}), (\ref{eq:uk}) and (\ref{bnd:nrj}) hold. We assume that $\uk$ is radially symmetrical for all $k\in\nn$. We assume that there exists $\delta_0\in (0,1)$ such that 
\bequa\label{hyp:lim:supuk}
\lim_{k\to +\infty}\sup_{B_{\delta_0}(0)}\uk=+\infty.
\eequa
Then for all $\delta\in (0,1)$ and for $k>0$ large enough, we have that $\sup_{B_\delta(0)}\uk=\uk(0)$.
\end{lem}

\smallskip\noindent{\it Proof of Lemma \ref{lem:uk:infty}:} It follows from (\ref{estim:wpe}) and (\ref{hyp:lim:supuk}) that for any $\delta\in (0,1)$, we have that $\lim_{k\to +\infty}\sup_{B_{\delta}(0)}\uk=+\infty$. It follows from the study of the monotonicity carried out in Step \ref{sec:prem}.1 that $\sup_{B_\delta(0)}\uk\in\{\uk(0),\uk(\delta)\}$. With (\ref{estim:wpe}), we get that there exists $C(\delta)>0$ such that $\uk(\delta)\leq C(\delta)$ for all $k\in\nn$. Since $\lim_{k\to +\infty}\sup_{B_{\delta}(0)}\uk=+\infty$, we get that the supremum is achieved at $0$ for $k>0$ large enough. This proves the Lemma \ref{lem:uk:infty}.
\hfill$\Box$

\medskip\noindent From now on, we assume that the sequence $(\uk)$ satisfies the hypothesis of Lemma \ref{lem:uk:infty}. In particular, we assume that for any $\delta\in (0,1)$, we have that
\bequa\label{max:uk}
\sup_{B_\delta(0)}\uk=\uk(0)\hbox{ and }\lim_{k\to +\infty}\uk(0)=+\infty.
\eequa

\medskip\noindent{\bf Step \ref{sec:prem}.5:} We now introduce the three fundamental types of convergence for $(E)$. This is a specificity of the bi-harmonic operator, compared to the Laplacian:

\begin{propdefi}\label{prop:defquad}
Let $(\Vk)_{k\in\nn}\in C^0(B)$ and $(\uk)_{k\in\nn}\in C^4(B)$ such that (\ref{hyp:Vk}), (\ref{eq:uk}) and (\ref{bnd:nrj}) hold. We assume that $\uk$ is radially symmetrical for all $k\in\nn$. We assume that (\ref{max:uk}) holds. We let
\bequa\label{def:vk}
\mk=e^{-\uk(0)}\hbox{ and }\vk(x):=\uk(\mk x)-\uk(0)
\eequa
for all $k\in\nn$ and all $x\in B_{\mk^{-1}}(0)$. Then one and only one of the following situations holds:\par
\smallskip(i: Log-convergence) For all $x\in\rn$,
$$\lim_{k\to +\infty}\vk(x)=\ln\frac{\sqrt{96}}{\sqrt{96}+|x|^2}.$$
Moreover, this convergence holds in $C^3_{loc}(\rn)$.\par

\smallskip(ii: Quadratic-Convergence (I).) There exists $a>0$, $v\in C^4(\rn)$ such that 
$$\Delta^2 v= e^{4v}\hbox{ in }\rn\hbox{ and }\lim_{|x|\to +\infty}\frac{v(x)}{|x|^2}=-a$$
and such that 
$$\lim_{k\to +\infty}\vk=v\hbox{ in }C^3_{loc}(\rn).$$

\smallskip(iii: Quadratic-Convergence (II).) We have that $\lim_{k\to +\infty}\Delta\vk(0)=+\infty$ and, for all $x\in\rn$,
$$\lim_{k\to +\infty}\frac{\vk(x)}{\Delta\vk(0)}\to -\frac{|x|^2}{8}.$$
Moreover, this convergence holds in $C^3_{loc}(\rn)$. 
\end{propdefi}

\medskip\noindent{\it Proof of Lemma \ref{prop:defquad}:} We let $\tVk(x):=\Vk(\mk x)$ for all $x\in B_{\mk^{-1}}(0)$ and all $k\in\nn$. In particular,
$$\lim_{k\to +\infty}\tVk=1\hbox{ in }C^0_{loc}(\rn).$$
Equation (\ref{eq:uk}) rewrites as
\bequa\label{eq:vk}
\Delta^2\vk=\tVk e^{4\vk}
\eequa
in $B_{\mk^{-1}}(0)$. Inequality (\ref{bnd:nrj}) rewrites as
\bequa\label{nrj:vk}
\int_{B_{\mk^{-1}}(0)}\tVk e^{4\vk}\, dx\leq\Lambda
\eequa
for all $k\in\nn$. Moreover, it follows from (\ref{max:uk}) and the definition (\ref{def:vk}) of $\vk$ that
\bequa\label{ineq:vk}
\vk(x)\leq \vk(0)=0
\eequa
for all $x\in B_{\mk^{-1}}(0)$. We let $R>0$. We proceed as in \cite{ars} and let $\wk\in C^4(B_R(0))$ such that
$$\left\{\begin{array}{ll}
\Delta^2\wk=\tVk e^{4\vk}&\hbox{ in }B_R(0)\\
\wk=\Delta\wk=0 &\hbox{ on }\partial B_R(0)
\end{array}\right\}.$$
It follows from standard elliptic theory that there exists $C(R)>0$ such that 
\bequa\label{est:wk:c3}
\Vert \wk\Vert_{C^{3,1/2}(\overline{B}_R(0))}\leq C(R)
\eequa
for all $k\in\nn$. We let $\varphi_k:=\vk-\wk$. It follows from (\ref{ineq:vk}) and (\ref{est:wk:c3}) that there exists $C(R)>0$ such that 
$$\varphi_k(x)\leq C(R)$$
for all $x\in B_R(0)$ and all $k\in\nn$. Proceeding as in \cite{ars}, we get that either $\varphi_k$ converges in $C^4_{loc}(B_R(0))$, or it converges in $C^4_{loc}(B_R(0))$ up to multiplication by a sequence. Coming back to the function $\vk=\wk+\varphi_k$ and using arbitrarily large $R>0$, we get that we are in one and only one of the following cases:





\smallskip\noindent{\it Case \ref{sec:prem}.5.1:} There exists $v\in C^3(\rn)$ such that 
$$\lim_{k\to +\infty}\vk=v\hbox{ in }C^3_{loc}(\rn).$$
Pasing to the limit in (\ref{eq:vk}), we get that $\Delta^2 v=e^{4v}$ in the distribution sense, and then $v\in C^4(\rn)$ by elliptic theory. With (\ref{ineq:vk}), we get that $v(x)\leq v(0)=0$ for all $x\in\rn$. Letting $k\to +\infty$ in (\ref{nrj:vk}), we get that $e^{4v}\in L^1(\rn)$. It follows from \cite{lin} (Theorem 1.1 and 1.2) that
$$\hbox{either }v(x)=\ln\frac{\sqrt{96}}{\sqrt{96}+|x|^2}\hbox{ or there exists }a>0\hbox{ such that }\lim_{|x|\to +\infty}\frac{v(x)}{|x|^2}=-a.$$
We recover (i) and (ii) of Proposition-Definition \ref{prop:defquad}. This ends Case \ref{sec:prem}.5.1.

\smallskip\noindent{\it Case \ref{sec:prem}.5.2:} There exists $\varphi\in C^3(\rn)$, there exists $(a_k)_{k\in\nn}\in\rr_{>0}$ such that $\lim_{k\to +\infty}a_k=+\infty$ and 
$$\lim_{k\to +\infty}\frac{\vk}{a_k}=\varphi\hbox{ in }C^3_{loc}(\rn).$$
Moreover, $\varphi\not\equiv 0$ and $\Delta^2\varphi=0$ in the distribution sense, and then $\varphi\in C^4(\rn)$ by elliptic theory. Passing to the limit in (\ref{ineq:vk}), we get that $\varphi(x)\leq \varphi(0)=0$ for all $x\in\rn$. It follows that there exists $\alpha>0$ such that $\varphi(x)=-\alpha|x|^2$ for all $x\in\rn$. Estimating $\Delta\vk(0)$, we get that
$$\lim_{k\to +\infty}\Delta\vk(0)=+\infty\hbox{ and }\lim_{k\to +\infty}\frac{\vk(x)}{\Delta\vk(0)}=-\frac{|x|^2}{8}$$
for all $x\in\rn$. Moreover, this convergence holds in $C^3_{loc}(\rn)$. In this case, we recover (iii) of Proposition-Definition \ref{prop:defquad}. This ends Case \ref{sec:prem}.5.2, and therefore the proof of Proposition-Definition \ref{prop:defquad}.

\hfill$\Box$

\medskip\noindent{\bf Step \ref{sec:prem}.6:} We state a very useful integral inequality. In the next section, this inequality will allow us to distinguish the three types of convergence above.

\begin{lem}\label{lem:intvk} 
Let $(\Vk)_{k\in\nn}\in C^0(B)$ and $(\uk)_{k\in\nn}\in C^4(B)$ such that (\ref{hyp:Vk}), (\ref{eq:uk}) and (\ref{bnd:nrj}) hold. We assume that $\uk$ is radially symmetrical for all $k\in\nn$. We assume that (\ref{max:uk}) holds. Then, for any $0<\delta<1$, there exists $C(\delta)>0$ such that 
\bequa\label{estim:intvk}
\int_{B_R(0)}|\Delta\vk-\mk^2\Delta\uk(\delta)|\, dx\leq C(\delta)R^2
\eequa
for all $k\in\nn$ and all $R<\delta\mk^{-1}$. In this expression, $\mk$ and $\vk$ are as in (\ref{def:vk}).
\end{lem}

\smallskip\noindent{\it Proof of Lemma \ref{lem:intvk}:} We follow the argument of Robert-Struwe \cite{robertstruwe}. We let $G_\delta$ be the Green's function for the Laplacian on $B_\delta(0)$ with Dirichlet boundary condition. Since $\Delta\uk$ is radially symmetrical, we get that
$$\Delta\uk(z)=\int_{B_\delta(0)}G_\delta(z,y)\Delta^2\uk(y)\, dy +\Delta\uk(\delta)$$
for all $z\in B_\delta(0)$. We choose $x\in \rn$ such that $|x|<\delta\mk^{-1}$. Using (\ref{def:vk}), we get that
\bequa\label{eq:useful}
\Delta\vk(x)-\mk^2\Delta\uk(\delta)= \int_{B_\delta(0)}\mk^2 G_\delta(\mk x,y)\Delta^2\uk(y)\, dy.
\eequa
Standard estimates on the Green's function (see for instance \cite{gr}) yield that there exists $C(\delta)>0$ such that
\bequa\label{est:G:d}
|G_\delta(x,y)|\leq \frac{C(\delta)}{|x-y|^2}
\eequa
for all $x,y\in B_{\delta}(0)$. Integrating (\ref{eq:useful}), using \eqref{est:G:d} and (\ref{eq:uk}), we get that
\beq
&&\int_{B_R(0)}|\Delta\vk-\mk^2\Delta\uk(\delta)|\, dx\\
&& \leq  \int_{x\in B_R(0)}\int_{y\in B_\delta(0)}\mk^2 G_\delta(\mk x,y)\Vk(y)e^{4\uk(y)}\, dy\, dx\\
&&\leq  C(\delta)\int_{B_\delta(0)} \Vk(y)e^{4\uk(y)}\left(\int_{B_R(0)}\frac{\mk^2}{|\mk x-y|^2}\, dx\right)\, dy\\
&&\leq  C(\delta)\int_{B_\delta(0)}\Vk e^{4\uk(y)}\left( CR^2\right)\, dy\leq  C(\delta) \Lambda  R^2,
\eeq
where $C(\delta)>0$ is independant of $k\in\nn$ and $R\in(0,\delta\mk^{-1})$. In this last inequality, we have used (\ref{bnd:nrj}). This proves Lemma \ref{lem:intvk}.
\hfill$\Box$

\medskip\noindent The key-quantity in Step \ref{sec:prem}.6 is the limit of $\mk^2\Delta\uk(\delta)$ when $k\to +\infty$. We separate the study in three cases, each of the following three sections is devoted to one of these cases. Thanks to them, we will recover the three notions of convergence of Proposition \ref{prop:defquad}.

\section{The case $\lim_{k\to +\infty}\mk^2\Delta\uk(\delta)=+\infty$}\label{sec:liminfty}

In this situation, we show that the second type of quadratic convergence of Proposition \ref{prop:defquad} holds and that $\Vk e^{4\uk}\rightharpoonup 0$ when $k\to +\infty$ in the sense of measures.

\medskip\noindent{\bf Step \ref{sec:liminfty}.1:} We prove that quadratic-convergence (II) of Proposition \ref{prop:defquad} holds in this case. More precisely,

\begin{lem}\label{lem:cvinfty} 
Let $(\Vk)_{k\in\nn}\in C^0(B)$ and $(\uk)_{k\in\nn}\in C^4(B)$ such that (\ref{hyp:Vk}), (\ref{eq:uk}) and (\ref{bnd:nrj}) hold. We assume that $\uk$ is radially symmetrical for all $k\in\nn$. We assume that (\ref{max:uk}) holds. We assume that there exists $\delta_0\in (0,1)$ such that
\bequa\label{hyp:sec:liminfty}
\lim_{k\to +\infty}\mk^2|\Delta\uk(\delta_0)|=+\infty,
\eequa
where $\mk$ is as in \eqref{def:vk}. Then, the second type of quadratic convergence of Proposition \ref{prop:defquad} holds. In addition, for any $\delta\in (0,1)$,
$$\lim_{k\to +\infty}\frac{\Delta\uk}{\Delta\uk(\delta)}= 1$$
in $C^1_{loc}(B)$ when $k\to +\infty$. 
\end{lem}

\smallskip\noindent{\it Proof of Lemma \ref{lem:cvinfty}:} Let $\delta_0\in (0,1)$ as in the Lemma. Let $R>0$. It follows from (\ref{estim:intvk}) and (\ref{hyp:sec:liminfty}) that
\bequa\label{eq:lim:case1}
\left\Vert\Delta\left(\frac{\vk}{\mk^2\Delta\uk(\delta_0)}\right)-1\right\Vert_{L^1(B_R(0))}=o(1)
\eequa
when $k\to +\infty$. It follows from (\ref{eq:vk}), (\ref{ineq:vk}) and (\ref{hyp:sec:liminfty}) that
$$\Delta\left(\Delta\left(\frac{\vk}{\mk^2\Delta\uk(\delta_0)}\right)-1\right)=\frac{\tVk e^{4\vk}}{\mk^2\Delta\uk(\delta_0)}=o(1)$$
where $o(1)\to 0$ in $C^0(B_R(0))$. It follows from \eqref{eq:lim:case1} and standard elliptic theory that
\bequa\label{eq:lim:vk}
\left\Vert\Delta\left(\frac{\vk}{\mk^2\Delta\uk(\delta_0)}\right)-1\right\Vert_{L^\infty(B_{R/2}(0))}\to 0
\eequa
when $k\to +\infty$. With (\ref{ineq:vk}), we get that there exists $\psi_R\in C^4(B_{R/2}(0))$ such that $\psi_R\not\equiv 0$ and 
$$\lim_{k\to +\infty}\frac{\vk}{\mk^2\Delta\uk(\delta_0)}=\psi_R\hbox{ in }C^3(B_{R/4}(0)).$$
Moreover, (\ref{ineq:vk}) yields $\psi_R(x)\leq\psi_R(0)=0$ for all $x\in B_{R/4}(0)$. With (\ref{eq:lim:vk}), we get that $\Delta\psi_R=1$. Since the functions are radial, we get that $\psi_R(x)=-|x|^2/8$ for all $x\in B_{R/4}(0)$. In particular, taking $R$ arbitrarily large, we get that
\bequa\label{cv:quad:2}
\lim_{k\to +\infty}\frac{\vk}{\mk^2\Delta\uk(\delta_0)}=-\frac{|x|^2}{8}\hbox{ in }C^3_{loc}(\rn).
\eequa
Computing the Laplacian of $\vk$ at $0$, we get that 
$$\lim_{k\to +\infty}\frac{\Delta\uk(0)}{\Delta\uk(\delta_0)}=1.$$
In particular, $\Delta\uk(\delta_0)>0$ for $k>0$ large and $\lim_{k\to +\infty}\Delta\uk(0)=+\infty$. Combining this limit with (\ref{cv:quad:2}), we obtain that the second type of quadratic convergence of Proposition \ref{prop:defquad} holds.

\medskip\noindent We let $\psi_k\in C^2(B)$ such that $\psi_k=\frac{\Delta\uk}{\Delta\uk(0)}$. With the equation (\ref{eq:uk}) and the estimate (\ref{estim:wpe}), we get that
$$\lim_{k\to +\infty}\Delta\psi_k=0\hbox{ in }C^0_{loc}(B\setminus\{0\}).$$
Since $\Delta\uk$ is decreasing, we have that $\psi_k(x)\leq 1$ for all $x\in B$. Noting that we have that $\lim_{k\to +\infty}\psi_k(\delta_0)=1$, it follows from elliptic theory that there exists $\psi\in C^2(B\setminus\{0\})$ such that $\lim_{k\to +\infty}\psi_k=\psi$ in $C^1_{loc}(B\setminus\{0\})$ and $\Delta\psi=0$. Letting $k\to +\infty$, we get that $\psi(\delta_0)=1$. Since $\psi\leq 1$ in $B$ and $\psi$ is non-increasing, we get that $\psi\equiv 1$. In addition, since $\psi_k$ is decreasing and achieves the value $1$ at $0$, we get that
$$\lim_{k\to +\infty}\frac{\Delta\uk}{\Delta\uk(0)}=1\hbox{ in }C^0_{loc}(B).$$
This ends the proof of the Lemma \ref{lem:cvinfty}.
\hfill$\Box$

\medskip\noindent{\bf Step \ref{sec:liminfty}.2:} In the case of quadratic convergence, the quadratic term happens to dominate the other ones asymptotically. More precisely, we have the following. Note that this Lemma does not use hypothesis \eqref{hyp:sec:liminfty}.

\begin{lem}\label{lem:EF1infty}[Pointwise estimate (I)] 
Let $(\Vk)_{k\in\nn}\in C^0(B)$ and $(\uk)_{k\in\nn}\in C^4(B)$ such that (\ref{hyp:Vk}), (\ref{eq:uk}) and (\ref{bnd:nrj}) hold. We assume that $\uk$ is radially symmetrical for all $k\in\nn$. We assume that (\ref{max:uk}) holds. Then for any $0<\delta<1$ there exists $C(\delta)>0$ such that
$$|x|\left|\nabla\uk(x)+\frac{\Delta\uk(\delta)}{4}x\right|\leq C(\delta)$$
for all $x\in B_{\delta}(0)$ and all $k\in\nn$.
\end{lem}

\smallskip\noindent{\it Proof of Lemma \ref{lem:EF1infty}:} We let $\delta\in (0,1)$. We let $H_\delta$ be the Green's function $\Delta^2$ on $B_\delta(0)$ with Navier condition, that is for any $x\in B_\delta(0)$,
$$\left\{\begin{array}{ll}
\Delta^2 H_\delta(x,\cdot)=\delta_x&\hbox{ in }{\mathcal D}'(B_\delta(0))\\
H_\delta(x,\cdot)=\Delta H_\delta(x,\cdot)=0 &\hbox{ on }\partial B_\delta(0)
\end{array}\right.$$
As easily checked, we have that $H_\delta=G_\delta*G_\delta$ where $*$ denotes the product of convolution and $G_\delta$ is the Green's function for $\Delta$ on $B_\delta(0)$ with Dirichlet boundary condition. Since $\uk$ is radially symmetrical, we get that
$$\uk(x)=\int_{B_\delta(0)}H_\delta(x,y)\Delta^2\uk(y)\, dy+\uk(\delta)+\frac{\delta^2-|x|^2}{8}\Delta\uk(\delta)$$
for all $x\in B_\delta(0)$. Differentiating this identity, we get that
\bequa\label{eq:temp:green}
\nabla\uk(x)=\int_{B_\delta(0)}\nabla H_\delta(x,y)\Delta^2\uk(y)\, dy-\frac{\Delta\uk(\delta)}{4}x
\eequa
for all $x\in B_{\delta}(0)$. Standard estimates on the Green's function (see for instance \cite{gr}) yield that there exists $C(\delta)>0$ such that
\bequa\label{est:green:1}
|\nabla H_\delta(x,y)|\leq \frac{C(\delta)}{|x-y|}
\eequa
for all $x,y\in B_{\delta}(0)$. Plugging (\ref{est:green:1}) into (\ref{eq:temp:green}), we get that
$$\left|\nabla\uk(x)+\frac{\Delta\uk(\delta)}{4}x\right|\leq C(\delta)\int_{B_\delta(0)}\frac{e^{4\uk(y)}}{|x-y|}\, dy$$
for all $x\in B_{\delta}(0)$. Using the pointwise estimate (\ref{estim:wpe}), we get that
\beq
\int_{B_\delta(0)}\frac{e^{4\uk(y)}}{|x-y|}\, dy&\leq & \int_{B_\delta(0)\cap B_{|x|/2}(0)}\frac{e^{4\uk(y)}}{|x-y|}\, dy+\int_{B_\delta(0)\setminus B_{|x|/2}(0)}\frac{e^{4\uk(y)}}{|x-y|}\, dy\\
&\leq & \frac{2}{|x|} \int_{B_\delta(0)}e^{4\uk(y)}\, dy+C\int_{B_\delta(0)\setminus B_{|x|/2}(0)}\frac{1}{|y|^4|x-y|}\, dy\\
&\leq & \frac{2\Lambda}{|x|}+\frac{C}{|x|}\int_{B_{\delta/|x|}(0)\setminus B_{1/2}(0)}\frac{1}{|y|^4\left|\frac{x}{|x|}-y\right|}\, dy\\
&\leq & \frac{C'(\delta,\Lambda)}{|x|}
\eeq
for all $x\in B_{\delta}(0)\setminus\{0\}$ and all $k\in\nn$. Here $C'(\delta,\Lambda)$ depends only on $\delta$ and $\Lambda$. This proves Lemma \ref{lem:EF1infty}.
\hfill$\Box$

\medskip\noindent{\bf Step \ref{sec:liminfty}.3:} We are in position to describe precisely the asymptotics of the $\uk$'s when $k\to +\infty$. This is the object of the following Lemma:

\begin{lem}\label{lem:EF2infty}[Pointwise estimate (II)] 
Let $(\Vk)_{k\in\nn}\in C^0(B)$ and $(\uk)_{k\in\nn}\in C^4(B)$ such that (\ref{hyp:Vk}), (\ref{eq:uk}) and (\ref{bnd:nrj}) hold. We assume that $\uk$ is radially symmetrical for all $k\in\nn$. We assume that (\ref{max:uk}) holds. We assume that there exists $\delta_0\in (0,1)$ such that
$$\lim_{k\to +\infty}\mk^2|\Delta\uk(\delta_0)|=+\infty.$$
We let $0<\delta<1$. Then
\bequa\label{asymp:1}
\uk\left(\frac{x}{\sqrt{\Delta\uk(\delta)}}\right)-\uk(0)=-\frac{|x|^2}{8}+O(1)\ln(2+|x|^2)
\eequa
for all $x\in B_{\delta\sqrt{\Delta\uk(0)}}(0)$ and all $k\in\nn$, where $O(1)$ denotes a function such that there exists $C(\delta)>0$ such that $|O(1)(x,k)|\leq C(\delta)$ for all $x\in B_{\delta\sqrt{\Delta\uk(0)}}(0)$ and all $k\in\nn$.
\end{lem}

\smallskip\noindent{\it Proof of Lemma \ref{lem:EF2infty}:} We let $x\in B_{\delta\sqrt{\Delta\uk(0)}}(0)$ such that $|x|>1$. We let $x_0=\frac{x}{|x|}$. With the pointwise estimate of Lemma \ref{lem:EF1infty}, we get that
\beqn
&&\uk\left(\frac{x}{\sqrt{\Delta\uk(\delta)}}\right)-\uk\left(\frac{x_0}{\sqrt{\Delta\uk(\delta)}}\right)\nonumber\\
&&= \int_0^1\frac{\partial}{\partial t}\left[\uk\left((1-t)\frac{x_0}{\sqrt{\Delta\uk(\delta)}}+t\frac{x}{\sqrt{\Delta\uk(\delta)}}\right)\right]\, dt\nonumber\\
&&= \frac{1}{\sqrt{\Delta\uk(\delta)}}\int_0^1(x-x_0)^i\partial_i\uk\left(\frac{(1-t)x_0+tx}{\sqrt{\Delta\uk(\delta)}}\right)\, dt\nonumber\\
&&= -\frac{1}{4}\int_0^1(x-x_0)^i((1-t)x_0+tx)_i\, dt\nonumber\\
&&+\frac{1}{\sqrt{\Delta\uk(\delta)}}\int_0^1(x-x_0)_i\left(\partial_i\uk\left(\frac{(1-t)x_0+tx}{\sqrt{\Delta\uk(\delta)}}\right)+\frac{\sqrt{\Delta\uk(\delta)}}{4}((1-t)x_0+tx)_i\right)\, dt\nonumber\\
&&=-\frac{|x|^2}{8}+\frac{|x_0|^2}{8}+O(1)\int_0^1\frac{|x-x_0|}{|(1-t)x_0+tx|}\, dt\nonumber\\
&&=-\frac{|x|^2}{8}+\frac{|x_0|^2}{8}+O(1)\int_0^1\frac{|x|-1}{t(|x|-1)+1}\, dt\nonumber\\
&&=-\frac{|x|^2}{8}+\frac{|x_0|^2}{8}+O(1)\ln|x|\label{cvquad:1}
\eeqn
where $O(1)$ is a function which is bounded with respect to both $x$ and $k\in\nn$. We claim that
\bequa\label{cvquad:2}
\lim_{k\to +\infty}\left(\uk\left(\frac{x}{\sqrt{\Delta\uk(\delta)}}\right)-\uk(0)\right)=-\frac{|x|^2}{8}
\eequa
for all $x\in\rn$, and that this convergence holds in $C^1_{loc}(\rn)$. We prove the claim. We write that 
$$\Delta\left(\uk\left(\frac{x}{\sqrt{\Delta\uk(\delta)}}\right)-\uk(0)\right)=\frac{\Delta\uk\left(\frac{x}{\sqrt{\Delta\uk(\delta)}}\right)}{\Delta\uk(\delta)},$$
for all $x\in B_\delta(0)$ and all $k\in\nn$. It follows from Lemma \ref{lem:cvinfty}, (\ref{max:uk}) and standard elliptic theory, that there exists $\varphi\in C^1(\rn)$ such that 
$$\lim_{k\to +\infty}\left(\uk\left(\frac{x}{\sqrt{\Delta\uk(\delta)}}\right)-\uk(0)\right)=\varphi(x)$$
for all $x\in\rn$ when $k\to +\infty$. Moreover, $\varphi\in C^2(\rn)$, $\Delta\varphi =1$ and $\varphi\leq\varphi(0)=0$. Since $\varphi$ is radially symmetrical, we get that $\varphi(x)=-\frac{|x|^2}{8}$. This proves the claim.

\medskip\noindent The asymptotic (\ref{asymp:1}) follows from (\ref{cvquad:1}) and (\ref{cvquad:2}). This proves Lemma \ref{lem:EF2infty}. 
\hfill$\Box$

\medskip\noindent{\bf Step \ref{sec:liminfty}.4:} We prove the vanishing of the $L^1-$norm of $e^{4\uk}$ when $k\to +\infty$.
\begin{lem}\label{lem:massinfty} 
Let $(\Vk)_{k\in\nn}\in C^0(B)$ and $(\uk)_{k\in\nn}\in C^4(B)$ such that (\ref{hyp:Vk}), (\ref{eq:uk}) and (\ref{bnd:nrj}) hold. We assume that $\uk$ is radially symmetrical for all $k\in\nn$. We assume that (\ref{max:uk}) holds. We assume that there exists $\delta_0\in (0,1)$ such that
$$\lim_{k\to +\infty}\mk^2|\Delta\uk(\delta_0)|=+\infty.$$
Then, for any $\delta\in (0,1)$, we have that
$$\int_{B_\delta(0)}e^{4\uk}\, dx\to 0$$
when $k\to +\infty$. In particular $\Vk e^{4\uk}\, dx\rightharpoonup 0$ when $k\to +\infty$ in the sense of measures.
\end{lem}

\smallskip\noindent{\it Proof of Lemma \ref{lem:massinfty}:} We let $\delta\in (0,1)$. With the definition (\ref{def:vk}) of $\mk$ and a change of variables, we get that
$$\int_{B_{\delta}(0)}e^{4\uk}\, dx=\frac{1}{\mk^4\Delta\uk(\delta)^2}\int_{B_{\delta\sqrt{\Delta\uk(\delta)}}(0)}e^{4\left(\uk\left(\frac{x}{\sqrt{\Delta\uk(\delta)}}\right)-\uk(0)\right)}\, dx.$$
Since $\lim_{k\to +\infty}\mk^2|\Delta\uk(\delta_0)|=+\infty$, Lemma \ref{lem:cvinfty} yields that $\lim_{k\to +\infty}\mk^2\Delta\uk(\delta)=+\infty$. The asymptotic (\ref{asymp:1}) of Lemma \ref{lem:EF2infty} yield the conclusion of the Lemma \ref{lem:massinfty}.
\hfill$\Box$

\medskip\noindent Point {\it (ii.c)} of Theorem \ref{thm:main} follows from Lemma \ref{lem:massinfty}.

\section{The case $\lim_{k\to +\infty}\mk^2|\Delta\uk(\delta)|=K_\delta>0$}\label{sec:limK}
In this situation, we show that the first type of quadratic convergence of Proposition \ref{prop:defquad} holds. Moreover, we describe the asymptotics for $\uk$.

\medskip\noindent{\bf Step \ref{sec:limK}.1:} We first prove that the quadratic-convergence (I) holds in this case. More precisely,
\begin{lem}\label{lem:cvK} 
Let $(\Vk)_{k\in\nn}\in C^0(B)$ and $(\uk)_{k\in\nn}\in C^4(B)$ such that (\ref{hyp:Vk}), (\ref{eq:uk}) and (\ref{bnd:nrj}) hold. We assume that $\uk$ is radially symmetrical for all $k\in\nn$. We assume that (\ref{max:uk}) holds. We assume that there exists $\delta_0\in (0,1)$ and $K_{\delta_0}\in\rr$ such that
$$\lim_{k\to +\infty}\mk^2|\Delta\uk(\delta_0)|=K_{\delta_0}>0.$$
Then the first type of quadratic convergence of Proposition-Definition \ref{prop:defquad} holds. In addition, we have that there exists $K>0$ such that
$$\lim_{k\to +\infty} \frac{\Delta\uk}{\Delta\uk(0)}= K\hbox{ in }C^1_{loc}(B\setminus\{0\}).$$
\end{lem}

\smallskip\noindent{\it Proof of Lemma \ref{lem:cvK}:} Let $R>0$. Since, up to a subsequence,
\bequa\label{lim:ueK}
\lim_{k\to +\infty}\mk^2\Delta\uk(\delta_0)=K_{\delta_0}\neq 0.
\eequa
It follows from (\ref{estim:intvk}) that
$$\left\Vert\Delta\vk\right\Vert_{L^1(B_R(0))}=O(1)$$
when $k\to +\infty$. It follows from equation (\ref{eq:vk}), inequation (\ref{ineq:vk}) and elliptic theory that
\bequa\label{eq:est:K}
\left\Vert\Delta\vk\right\Vert_{C^1(B_{R/2}(0))}=O(1)
\eequa
when $k\to +\infty$. Inequation (\ref{ineq:vk}), equations \eqref{eq:vk} and \eqref{eq:est:K}, the Harnack inequality and standard elliptic theory yield that there exists $v\in C^3(\rn)$ such that
\bequa\label{cvK:2}
\vk\to v\hbox{ in }C^3_{loc}(\rn),
\eequa
where $\Delta^2 v=e^{4v}$ in the distribution sense. Elliptic theory yields that $v\in C^4(\rn)$. We are then in Case (i) or (ii) of Proposition-Definition \ref{prop:defquad}. 

\medskip\noindent We claim that we are in Case (ii) of Proposition-Definition \ref{prop:defquad}. We proceed by contradiction and assume that Case (i) of Proposition-Definition \ref{prop:defquad} holds. We then get that $v=v_0$ where 
$$v_0(x)=\ln\frac{\sqrt{96}}{\sqrt{96}+|x|^2}$$
for all $x\in\rn$. We let $R>0$. We let $k\to +\infty$ in (\ref{estim:intvk}) and get that 
$$\int_{B_R(0)}|\Delta v(x)-K_{\delta_0}|\, dx\leq C R^2.$$
Since $v=v_0$, using the explicit expression of $v_0$ above and letting $R\to +\infty$, we get that there exists a constant $C>0$ independant of $R>0$ such that
$$|K_{\delta_0}|\leq C R^{-2}$$
for all $R>0$. Letting $R\to +\infty$, we get that $K_{\delta_0}=0$. A contradiction with our initial assumption \eqref{lim:ueK}. Then Case (i) does not hold and we are in Case (ii). 

\medskip\noindent It follows from Case (ii) of Proposition-Definition \ref{prop:defquad} and Theorem 1.2 of \cite{lin} that there exists $a>0$ such that
\bequa\label{asymp:deltav}
\lim_{|x|\to +\infty}\frac{v(x)}{|x|^2}=-a\hbox{ and }\lim_{|x|\to +\infty}\Delta v(x)=8a.
\eequa
We let $\delta\in (0,1)$. With (\ref{estim:intvk}), we get that there exists $C(\delta)>0$ such that
\bequa\label{eq:proof:quad2}
\int_{B_R(0)}|\Delta\vk-\mk^2\Delta\uk(\delta)|\, dx\leq C(\delta)R^2
\eequa
for all $R\in (0,\delta\mk^{-1})$. It the follows from (\ref{cvK:2}) that there exists $K_\delta\in\rr$ such that $\lim_{k\to +\infty}\mk^2\Delta\uk(\delta)=K_\delta$. Passing to the limit $k\to +\infty$ in (\ref{eq:proof:quad2}), we get that
$$\int_{B_R(0)}|\Delta v-K_\delta|\, dx\leq C(\delta)R^2
$$
for all $R>0$. Letting $R\to +\infty$ in this inequality and using (\ref{asymp:deltav}), we get that $K_\delta=8a>0$ for all $\delta\in (0,1)$. In particular, with (\ref{lim:ueK}) and (\ref{cvK:2}), we get that there exists $K>0$ such that for any $\delta\in (0,1)$,
$$\lim_{k\to +\infty}\frac{\Delta\uk(\delta)}{\Delta\uk(0)}=K>0.$$
The last assertion of Lemma \ref{lem:cvK} follows from this limit, equation $(E)$, inequality (\ref{estim:wpe}) the decreasing of $\Delta\uk$ and standard elliptic theory.
\hfill$\Box$

\medskip\noindent{\bf Step \ref{sec:limK}.2:} With some arguments very similar to the ones developed in the proof of Lemma \ref{lem:EF2infty}, we get the following Lemma. We omit the proof:

\begin{lem}\label{lem:EFK} 
Let $(\Vk)_{k\in\nn}\in C^0(B)$ and $(\uk)_{k\in\nn}\in C^4(B)$ such that (\ref{hyp:Vk}), (\ref{eq:uk}) and (\ref{bnd:nrj}) hold. We assume that $\uk$ is radially symmetrical for all $k\in\nn$. We assume that (\ref{max:uk}) holds. We assume that there exists $\delta_0\in (0,1)$ such that
$$\lim_{k\to +\infty}\mk^2|\Delta\uk(\delta_0)|=K_{\delta_0}>0.$$
We let $0<\delta<1$. Then there exists a sequence $(a_k)_{k\in\nn}\in\rr$ such that $\lim_{k\to +\infty}a_k=a_\infty>0$ and such that
$$\vk(x)=-a_k|x|^2+O(1)\ln(2+|x|^2)$$
for all $x\in B_{\delta\mk^{-1}}(0)$ and all $k\in\nn$, where $O(1)$ denotes a function such that there exists $C(\delta)>0$ such that $|O(1)(x,k)|\leq C(\delta)$ for all $x\in B_{\delta\mk^{-1}}(0)$ and all $k\in\nn$.
\end{lem}

\medskip\noindent As a consequence of this pointwise estimate, we get the following quantization of the $L^1-$norm of $e^{4\uk}$:

\begin{lem}\label{lem:massK} 
Let $(\Vk)_{k\in\nn}\in C^0(B)$ and $(\uk)_{k\in\nn}\in C^4(B)$ such that (\ref{hyp:Vk}), (\ref{eq:uk}) and (\ref{bnd:nrj}) hold. We assume that $\uk$ is radially symmetrical for all $k\in\nn$. We assume that (\ref{max:uk}) holds. We assume that there exists $\delta_0\in (0,1)$ such that
$$\lim_{k\to +\infty}\mk^2|\Delta\uk(\delta_0)|=K_{\delta_0}>0.$$
Then for any $\delta\in (0,1)$, we have that
$$\lim_{k\to +\infty}\int_{B_\delta(0)}\Vk  e^{4\uk}\, dx=\int_{\rn}e^{4v}\, dx<16\pi^2.$$
In other words, $\Vk e^{4\uk}\, dx\rightharpoonup (\int_{\rn}e^{4v}\, dx)\delta_0$ when $k\to +\infty$ in the sense of the measures.
\end{lem}

\smallskip\noindent{\it Proof of Lemma \ref{lem:massK}:} It follows from Lemma \ref{lem:EFK} that there exists $C=C(\delta)>0$ such that
$$\vk(x)\leq -\frac{a_\infty}{2}|x|^2+C$$
for all $x\in B_{\delta\mk^{-1}}(0)$ and all $k\in\nn$. We let $R>0$. With a change of variable, we get that
$$\int_{B_\delta(0)\setminus B_{R\mk}(0)}\Vk e^{4\uk}\, dx=\int_{B_{\delta/{\mk}}(0)\setminus B_R(0)}\tVk e^{4\vk}\, dx\leq 2\int_{\rn\setminus B_R(0)}e^{-2a_\infty|x|^2+4C}\, dx.$$
As a consequence,
\bequa\label{eq:lim:K:1}
\lim_{R\to +\infty}\lim_{k\to +\infty}\int_{B_\delta(0)\setminus B_{R\mk}(0)}\Vk e^{4\uk}\, dx=0.
\eequa
On the other hand, with a change of variables and letting $k\to +\infty$, we get that
\bequa\label{eq:lim:K:2}
\int_{ B_{R\mk}(0)}\Vk e^{4\uk}\, dx=\int_{ B_{R}(0)}e^{4\vk}\, dx=\int_{ B_{R}(0)}e^{4v}\, dx+o(1)
\eequa
when $k\to +\infty$. Summing \eqref{eq:lim:K:1} and \eqref{eq:lim:K:2} and letting $k\to +\infty$ and then $R\to +\infty$, we get that
$$\lim_{k\to +\infty}\int_{B_\delta(0)}\Vk  e^{4\uk}\, dx=\int_{\rn}e^{4v}\, dx.$$
Moreover, it follows from \cite{lin}, Theorem 1.2, that
$$\int_{\rn}e^{4v}\, dx<16\pi^2.$$
This ends the proof of Lemma \ref{lem:massK}.
\hfill$\Box$

\medskip\noindent Point {\it (ii.b)} of Theorem \ref{thm:main} follows from Lemma \ref{lem:massK}.




\section{The case $\lim_{k\to +\infty}\mk^2\Delta\uk(\delta)=0$}\label{sec:lim0}
In this case, the behaviour of the $\uk$'s is much more standard and is similar to the two-dimensional corresponding problem. We show that the Log-convergence of Proposition-Definition \ref{prop:defquad} holds. Moreover, we describe the asymptotics for $\uk$.

\medskip\noindent{\bf Step \ref{sec:lim0}.1:} We first prove that the Log-convergence holds in this case. More precisely,
\begin{lem}\label{lem:cv0} 
Let $(\Vk)_{k\in\nn}\in C^0(B)$ and $(\uk)_{k\in\nn}\in C^4(B)$ such that (\ref{hyp:Vk}), (\ref{eq:uk}) and (\ref{bnd:nrj}) hold. We assume that $\uk$ is radially symmetrical for all $k\in\nn$. We assume that (\ref{max:uk}) holds. We assume that there exists $\delta_0\in (0,1)$ such that
$$\lim_{k\to +\infty}\mk^2\Delta\uk(\delta_0)=0.$$
Then for any $x\in\rn$,
$$\lim_{k\to +\infty}\vk(x)= \ln\frac{\sqrt{96}}{\sqrt{96}+|x|^2},$$
where $\vk$ is as in \eqref{def:vk}. Moreover, this convergence holds in $C^3_{loc}(\rn)$.
\end{lem}
\smallskip\noindent{\it Proof of Lemma \ref{lem:cv0}:} With some arguments similar to the ones developed in the proof of Lemma \ref{lem:cvK}, we get that there exists $v\in C^4(\rn)$ such that $\lim_{k\to +\infty}\vk=v$ in $C^3_{loc}(\rn)$. Moreover, $\Delta^2 v= e^{4v}$ and $e^v\in L^1(\rn)$. We are then in Case (i) or (ii) of  Proposition-Definition \ref{prop:defquad}. We let $k\to +\infty$ in (\ref{estim:intvk}) and get for any $R>0$
in $\rn$ that
\bequa\label{est:log:v}
\int_{B_R(0)}|\Delta v(x)|\, dx\leq C R^2.
\eequa
We assume by contradiction that Case (ii) holds. It then follows from Lin \cite{lin} that $\lim_{|x|\to +\infty}\Delta v(x)=8a>0$. Letting $R\to +\infty$ in (\ref{est:log:v}), we get that $8a=0$. A contradiction. We are then in Case (i) of Proposition-Definition \ref{prop:defquad} and $v(x)=\ln\frac{\sqrt{96}}{\sqrt{96}+|x|^2}$ for all $x\in\rn$, that is Log-Convergence holds. This proves Lemma \ref{lem:cv0}.
\hfill$\Box$

\medskip\noindent A consequence of this Lemma is the following. With a change of variable, we get that
$$\int_{ B_{R\mk}(0)}\Vk e^{4\uk}\, dx=\int_{ B_{R}(0)}\tVk e^{4\vk}\, dx=\int_{ B_{R}(0)}e^{4v}\, dx+o(1)$$
when $k\to +\infty$. Passing to the limit $k\to +\infty$ and then $R\to+\infty$, we get that
\bequa\label{lim:mass0}
\lim_{R\to +\infty}\lim_{k\to +\infty}\int_{ B_{R\mk}(0)}\Vk e^{4\uk}\, dx=16\pi^2.
\eequa

\medskip\noindent{\bf Step \ref{sec:lim0}.2:} We are in position to deal with the convergence outside $0$. This is the object of the following Lemma:

\begin{lem}\label{lem:infty} 
Let $(\Vk)_{k\in\nn}\in C^0(B)$ and $(\uk)_{k\in\nn}\in C^4(B)$ such that (\ref{hyp:Vk}), (\ref{eq:uk}) and (\ref{bnd:nrj}) hold. We assume that $\uk$ is radially symmetrical for all $k\in\nn$. We assume that (\ref{max:uk}) holds. We assume that there exists $\delta_0\in (0,1)$ such that
$$\lim_{k\to +\infty}\mk^2\Delta\uk(\delta_0)=0.$$
Then $\lim_{k\to +\infty}\uk= -\infty$ uniformly on every compact subset of $B\setminus\{0\}$. 
\end{lem}

\smallskip\noindent{\it Proof of Lemma \ref{lem:infty}:} Assume that the conclusion is false. It then follows from Lemma \ref{cv:out0} that for any $K\subset\subset B\setminus\{0\}$, there exists $C(K)>0$ such that 
\bequa\label{est:uk:out}
|\uk(z)|+|\Delta\uk(z)|\leq C(K)
\eequa
for all $z\in K$. We let $\delta\in (0,1/2)$ and we let $H_\delta$ be the Green's function for $\Delta^2$ in $B_{\delta}(0)$ with Navier condition on the boundary, that is for any $x\in B_\delta(0)$, we have that
$$\Delta^2 H_\delta(x,\cdot)=\delta_x$$
for all $x\in B_{\delta}(x)$ and $H_\delta(x,\cdot)=\Delta H_{\delta}(x,\cdot)=0$ on $\partial B_{\delta}(x)$. We let $x\in B_{\delta}(0)\setminus\{0\}$. Since $\uk$ is radially symmetrical, we have that
$$\uk(x)=\int_{B_\delta(0)}H_\delta(x,y)\Vk(y) e^{4\uk(y)}\, dy+\uk(\delta)+\frac{\delta^2-|x|^2}{8}\Delta\uk(\delta).$$ 
We let $\alpha>0$ small. Since $\uk$ is uniformly bounded in $L^\infty$ outside $0$ and since $H_\delta>0$, we get with (\ref{hyp:Vk}), (\ref{lim:mass0}) and \eqref{est:uk:out} that there exists $C>0$ independant of $x$ and $\alpha>0$ such that
\beq
\uk(x) & \geq  & \int_{B_{R\mk}(0)}H_\delta(x,y) \Vk(y) e^{4\uk(y)}\, dy-C\\
&\geq & \int_{B_{R}(0)}H_\delta(x,\mk y) \tVk(y)e^{4\vk(y)}\, dy-C\\
&\geq & \int_{B_{R}(0)}H_\delta(x,0)\lim_{k\to +\infty}\left(\tVk(y) e^{4\vk(y)}\right)\, dy-C+o(1)\\
&\geq & 16\pi^2 H_\delta(x,0)-C+o(1)
\eeq
for $x\in B$ such that $|x|\geq\alpha$ and for $k$ large enough depending only on $\alpha$. Since $H_\delta(x,0)=\frac{1}{8\pi^2}\ln\frac{\delta}{|x|}+\frac{|x|^2-\delta^2}{32\pi^2\delta^2}$ for $x\in B_{\delta/2}(x_0)$. We then get that
$$\uk(x)\geq 2\ln\frac{1}{|x|}-C'+o(1)$$
for $x\in B_{\delta}(0)\setminus B_{\alpha}(0)$ and $k$ large depending only on $\alpha>0$. We then get that for any $0<\alpha<\beta$ small,
$$\Lambda\geq \int_{B_\beta(0)\setminus B_\alpha(0)} \Vk e^{4\uk}\, dx\geq C\int_{B_\beta(0)\setminus B_\alpha(0)}\frac{1}{|x|^8}\, dx.$$
We get a contradiction by letting $\alpha\to 0$. Then $\uk\to -\infty$ on compact subsets of $B\setminus\{0\}$ when $k\to +\infty$ and Lemma \ref{lem:infty} is proved.
\hfill$\Box$

\medskip\noindent{\bf Step \ref{sec:lim0}.3:} We now prove that the whole $L^1-$norm of $e^{4\uk}$ is actually $16\pi^2$. We borrow ideas from Schoen-Zhang \cite{schoenzhang}, Druet \cite{druet} and Druet-Robert \cite{druetrobert}.

\begin{lem}\label{lem:mono1} 
Let $(\Vk)_{k\in\nn}\in C^0(B)$ and $(\uk)_{k\in\nn}\in C^4(B)$ such that (\ref{hyp:Vk}), (\ref{eq:uk}) and (\ref{bnd:nrj}) hold. We assume that $\uk$ is radially symmetrical for all $k\in\nn$. We assume that (\ref{max:uk}) holds. We assume that there exists $\delta_0\in (0,1)$ such that
$$\lim_{k\to +\infty}\mk^2\Delta\uk(\delta_0)=0.$$
We let $\delta\in (0,1)$. Then there exists $(r_k)_{k\in\nn}\in \rr_{>0}$ such that $r_k\in [0,\delta]$ for all $k\in\nn$ and

\smallskip(i) $\lim_{k\to +\infty}\frac{\rk}{\mk}=+\infty$,\par
(ii) $r\mapsto r e^{\uk(r)}$ is decreasing on $[4\mk,\rk]$,\par
(iii) $\uk\to -\infty$ uniformly on $\overline{B}_\delta(0)\setminus B_{\rk}(0)$.
\end{lem}

\smallskip\noindent{\it Proof of Lemma \ref{lem:mono1}:} We let $\delta\in (0,1)$. Without loss of generality, we assume that $\lim_{k\to +\infty}\mk^2\Delta\uk(\delta)=0$ (otherwise, we are back to the previous cases).

\medskip\noindent{\it Step \ref{sec:lim0}.3.1:} We claim that for any $R>4$, we have that
$$r\mapsto r e^{\uk(r)}\hbox{ is decreasing on }[4\mk,R\mk]$$
for $k$ large enough. Indeed, we let $r\in [4\mk,R\mk]$ and we let $\rho_k:=\frac{r}{\mk}$. With Lemma \ref{lem:cv0}, we have that
\beq
(r e^{\uk(r)})'(r)&=&\mk^{-1}\frac{d}{dr}\left(r\mk e^{\uk(r\mk)}\right)(\rho_k)= \mk^{-1}\frac{d}{dr}\left(r e^{\vk(r)}\right)(\rho_k)\\
&=& \mk^{-1}\left(\frac{d}{dr}\left(r e^{v(r)}\right)+o(1)\right)(\rho_k)=\frac{\sqrt{96}}{\mk}\left(\frac{\sqrt{96}-\rho_k^2}{(\sqrt{96}+\rho_k^2)^2}+o(1)\right)
\eeq
where $o(1)\to 0$ when $k\to +\infty$ uniformly for $r\in [4\mk,R\mk]$. Since $\rho_k\geq 4$, the RHS in negative. Then $(r e^{\uk(r)})'<0$ and the function $r\mapsto r e^{\uk(r)}$ is decreasing on $[4\mk,R\mk]$.

\medskip\noindent{\it Step \ref{sec:lim0}.3.2:} We assume that $r\to  re^{\uk(r)}$ is decreasing on $[4\mk,\delta]$ for all $k\in\nn$. Then the conclusion of the Lemma holds with $\rk:=\delta$, and Lemma \ref{lem:mono1} is proved.

\medskip\noindent From now on, we assume that 
\bequa\label{hyp:non:dec}
r\to r e^{\uk(r)}\hbox{ is not decreasing on }[4\mk,\delta].
\eequa
We let
$$\rk:=\inf\{\rho\in [4\mk,\delta]/\, (r e^{\uk(r)})'(\rho)=0.\}$$

\medskip\noindent{\it Step \ref{sec:lim0}.3.3:} We claim that
\bequa\label{def:re}
\lim_{k\to +\infty}\frac{\rk}{\mk}=+\infty,\quad (r e^{\uk(r)})'(r)<0\hbox{ when }4\mk<r<\rk\hbox{ and }(r e^{\uk(r)})'(\rk)=0.
\eequa
Indeed, it follows from Step \ref{sec:lim0}.2.1 and (\ref{hyp:non:dec}) that $\rk$ is defined and satisfies the two last statements of (\ref{def:re}). The first statement is a consequence of Step \ref{sec:lim0}.2.1.

\medskip\noindent{\it Step \ref{sec:lim0}.3.4:} We claim that
\bequa\label{lim:expre}
\lim_{k\to +\infty}\rk e^{\uk(\rk)}=0.
\eequa
Indeed, we let $R\geq 4$. It follows from (\ref{def:re}) that $r e^{\uk(r)}$ is decreasing on $[R\mk,\rk]$. We then get that 
\beq
\rk e^{\uk(\rk)}&\leq & R\mk e^{\uk(R\mk)}\leq  R e^{\vk(R)}\\
&\leq & \left(R e^{v(R)}+o(1)\right)\leq \left(\frac{R\sqrt{96}}{\sqrt{96}+R^2}+o(1)\right)
\eeq
where $o(1)\to 0$ when $k\to +\infty$. Letting $k\to +\infty$ and then $R\to +\infty$, we get (\ref{lim:expre}). This ends Step \ref{sec:lim0}.2.4.

\medskip\noindent We let 
\bequa\label{def:tue}
\tuk(x)=\uk(\rk x)-\uk(\rk)
\eequa
for all $k\in\nn$ and all $x\in B_{\rk^{-1}}(0)$. 

\medskip\noindent{\it Step \ref{sec:lim0}.3.5:} We claim that there exists $a\geq 1$ such that for any $x\in\rn\setminus\{0\}$, we have that
\bequa\label{lim:tue2}
\lim_{k\to +\infty}\uk(\rk x)-\uk(\rk)=a\ln\frac{1}{|x|}+\frac{a-1}{2}(|x|^2-1).
\eequa
Moreover, this convergence holds in $C^3_{loc}(\rn\setminus\{0\})$. Indeed, equation (\ref{eq:uk}) rewrites as
\bequa\label{eq:tue}
\Delta^2\tuk(x)=\Vk(\rk x) \rk^4 e^{4\uk(\rk x)}=\Vk(\rk x) \rk^4 e^{4\uk(\rk)}e^{\tuk(x)}
\eequa
for all $k\in\nn$ and all $x\in B_{\rk^{-1}}(0)$. The system (\ref{def:re}) yields that
\bequa\label{hyp:tue}
r\tuk'(r)\leq -1\hbox{ for }\frac{4\mk}{\rk}\leq r\leq 1\hbox{ and }\tuk'(1)=-1.
\eequa
Proceeding as in Lemma \ref{lem:EF1infty} and using the definition (\ref{def:tue}), we get that there exists $C=C(\delta)>0$ such that
\bequa\label{estim:tue}
|\tuk'(r)+\frac{\Delta\uk(\delta)}{4}\rk^2 r|\leq \frac{C}{r}
\eequa
for all $k\in\nn$ and all $r\in (0,\delta\rk^{-1})$. Taking $r=1$ in (\ref{estim:tue}) and using (\ref{hyp:tue}), we then obtain that, up to a subsequence, there exists $\rho\in\rr$ such that
\bequa\label{lim:rho}
\lim_{k\to +\infty}\rk^2\Delta\uk(\delta)=\rho.
\eequa
Since $\tuk(1)=0$, it follows from (\ref{estim:tue}) and (\ref{lim:rho}) that for any $U\subset\subset \rn\setminus\{0\}$, there exists $C'(U)>0$ such that
$$|\tuk(x)|\leq C'(U)$$ 
for all $x\in U$ and all $k\in\nn$. It then follows from (\ref{eq:tue}), (\ref{lim:expre}) and standard elliptic theory that there exists $\tu\in C^4(\rn\setminus\{0\})$ such that $\Delta^2\tu=0$ and 
\bequa\label{lim:tue1}
\tuk\to\tu
\eequa
in $C^3_{loc}(\rn\setminus\{0\})$ when $k\to +\infty$. Since $\tu$ is radially symmetrical, we get that there exist $a,b,c,d\in\rr$ such that
\bequa\label{eq:tu1}
\tu(x)=a\ln\frac{1}{|x|}+\frac{b}{|x|^2}+c |x|^2+d
\eequa
for all $x\in \rn\setminus\{0\}$. Passing to the limit in (\ref{estim:tue}) and using \eqref{lim:rho}, we get that
\bequa\label{estim:tu}
|\tu'(r)+\frac{\rho}{4}r|\leq \frac{C}{r}
\eequa
for all $r>0$. It follows from (\ref{eq:tu1}) and (\ref{estim:tu}) that $b=2c+\frac{\rho}{4}=0$, so that we can write
$$\tu(x)=a\ln\frac{1}{|x|}-\frac{\rho}{8} |x|^2+d$$
for all $x\in \rn\setminus\{0\}$. Passing to the limit in (\ref{hyp:tue}), we get that
$$r\tu'(r)\leq -1\hbox{ for }r<1\hbox{ and }\tu'(1)=-1.$$
With the explicit expression (\ref{eq:tu1}) of $\tu$, we get that
$$a-1=\frac{|\rho|}{4}\geq 0\hbox{ and }\rho\leq 0.$$
Since $\tu(1)=1$, the claim follows.

\medskip\noindent{\it Step \ref{sec:lim0}.3.6:} We claim that 
$$a\geq 2.$$
Indeed, integrating by parts, we get that
\beqn\label{poho:1}
&&\int_{B_{\rk}(0)}x^i\partial_i\uk\Delta^2\uk\, dx\nonumber\\
&&=\int_{\partial B_{\rk}(0)}\left(\frac{(x,\nu)}{2}(\Delta\uk)^2+\Delta\uk\frac{\partial (x,\nabla\uk)}{\partial \nu}-(x,\nabla\uk)\frac{\partial \Delta\uk}{\partial\nu}\right)\, d\sigma
\eeqn
where $\nu$ denotes the outer normal vector at $\partial B_{\rk}(0)$. Using the change of variable $y=\rk x$ and the convergence (\ref{lim:tue2}), we get that
\beqn
&&\int_{\partial B_{\rk}(0)}\left(\frac{(x,\nu)}{2}(\Delta\uk)^2+\Delta\uk\frac{\partial (x,\nabla\uk)}{\partial \nu}-(x,\nabla\uk)\frac{\partial \Delta\uk}{\partial\nu}\right)\, d\sigma\nonumber\\
&&= \int_{\partial B_{1}(0)}\left(\frac{(x,\nu)}{2}(\Delta\tuk)^2+\Delta\tuk\frac{\partial (x,\nabla\tuk)}{\partial \nu}-(x,\nabla\tuk)\frac{\partial \Delta\tuk}{\partial\nu}\right)\, d\sigma\nonumber\\
&&=\int_{\partial B_{1}(0)}\left(\frac{(x,\nu)}{2}(\Delta\tu)^2+\Delta\tu\frac{\partial (x,\nabla\tu)}{\partial \nu}-(x,\nabla\tu)\frac{\partial \Delta\tu}{\partial\nu}\right)\, d\sigma+o(1)\nonumber\\
&&=-4\pi^2a^2+o(1)\label{poho:2}
\eeqn
where $o(1)\to 0$ when $k\to +\infty$. On the other hand, using (\ref{eq:uk}), we have that
\bequa\label{ipp:poho:A}
\int_{B_{\rk}(0)}x^i\partial_i\uk\Delta^2\uk\, dx=\int_{B_{\rk}(0)}x^i\partial_i\uk e^{4\uk}\, dx+\int_{B_{\rk}(0)}(\Vk-1)x^i\partial_i\uk e^{4\uk}\, dx.
\eequa
It follows from Lemma \ref{lem:EF1infty} and (\ref{lim:rho}) that there exists $C>0$ such that
\bequa\label{est:poho}
|x^i\partial_i\uk(x)|\leq C
\eequa
for all $x\in B_{\rk}(0)$. Properties (\ref{hyp:Vk}), (\ref{bnd:nrj}) and (\ref{est:poho}) yield 
\bequa\label{ipp:poho:B}
\lim_{k\to +\infty}\int_{B_{\rk}(0)}(\Vk-1)x^i\partial_i\uk e^{4\uk}\, dx=0.
\eequa
Plugging (\ref{poho:2}) and (\ref{ipp:poho:B}) into \eqref{poho:1} and (\ref{ipp:poho:A}), we get that
$$\int_{B_{\rk}(0)}x^i\partial_i\uk e^{4\uk}\, dx=-4\pi^2a^2+o(1)$$
when $k\to +\infty$. Integrating by parts, we get that
\beqn
-4\pi^2a^2&&=\int_{B_{\rk}(0)}x^i\partial_i\frac{e^{4\uk}}{4}\, dx+o(1)\nonumber\\
&&=-\int_{B_{\rk}(0)} e^{4\uk}\, dx+\int_{\partial B_{\rk}(0)}\frac{(x,\nu)}{4}e^{4\uk}\, d\sigma+o(1)\nonumber\\
&&=-\int_{B_{\rk}(0)} e^{4\uk}\, dx+\rk^4 e^{4\uk(\rk)}\int_{\partial B_{1}(0)}\frac{(x,\nu)}{4}e^{4\tuk}\, d\sigma+o(1).
\eeqn
With (\ref{lim:expre}) and (\ref{lim:tue1}), we then get that
$$\int_{B_{\rk}(0)} e^{4\uk}\, dx=4\pi^2a^2+o(1)$$
where $o(1)\to 0$ when $k\to +\infty$. Since (\ref{lim:mass0}) and (\ref{def:re}) hold, we then get that $|a|^2\geq 4$. Since $a\geq 1$, we get that $a\geq 2$, and the claim is proved.

\medskip\noindent{\it Step \ref{sec:lim0}.3.7:} We let $\delta\in (0,1)$. We claim that
$$\lim_{k\to +\infty}\sup_{[\rk,\delta]}\uk= -\infty.$$
Indeed, it follows from (\ref{lim:tue2}) and (\ref{lim:tue1}) that $r\tu'(r)=-a+(a-1)r^2$. Since $a>1$, we get that $\tu$ is decreasing on $\left(0,\sqrt{\frac{a}{a-1}}\right)$ and increasing on $\left(\sqrt{\frac{a}{a-1}},+\infty\right)$. It follows from the study of the monotonicity of $\uk$ provided in Case \ref{sec:prem}.2.3 of Step \ref{sec:prem}.2 that there exists $\tau_k\in (0,\delta)$ such that $\uk$ decreases on $(0,\tau_k)$ and increases on $(\tau_k,\delta)$. Since (\ref{lim:tue2}) and (\ref{lim:tue1}) hold and since the monotonicity of $\tu$ changes at $\sqrt{\frac{a}{a-1}}$, we get that 
\bequa\label{lim:tk}
\lim_{k\to +\infty}\frac{\tau_k}{\rk}=\sqrt{\frac{a}{a-1}}.
\eequa
We let $y_k\in \overline{B}_{\delta}(0)\setminus B_{\rk}(0)$ such that
$$\sup_{\overline{B}_{\delta}(0)\setminus B_{\rk}(0)}\uk=\uk(\yk).$$
We distinguish two cases:

\smallskip\noindent Case \ref{sec:lim0}.3.7.1: we assume that $\lim_{k\to +\infty}\frac{|\yk|}{\rk}=+\infty$. Then with \eqref{lim:tk}, we get that $\uk$ increases on $[\tau_k,\delta]$, and then $\uk(\yk)\leq \uk(\delta)$. With Lemma \ref{lem:infty}, we then get that $\lim_{k\to +\infty}\uk(\yk)=-\infty$.

\smallskip\noindent Case \ref{sec:lim0}.3.7.2: we assume that $|\yk|=O(\rk)$ when $k\to +\infty$. We let $z_k=\frac{\yk}{\rk}$. Since $|\yk|\geq \rk$, we get that, up to a subsequence, $\lim_{k\to +\infty}z_k=z_\infty\neq 0$. With (\ref{def:tue}), (\ref{lim:tue1}) and Case \ref{sec:lim0}.3.7.1, we get that
\beq
\uk(\yk)&=&\uk(\yk)-\uk(\tau_k)+\uk(\tau_k)\\
&\leq&\tuk(z_k)-\tuk\left(\frac{\tau_k}{\rk}\right)+\uk(\tau_k)\leq O(1)+\uk(\delta)
\eeq
and then with Lemma \ref{lem:infty}, we get that $\lim_{k\to +\infty}\uk(\yk)=-\infty$. This proves the claim.

\smallskip\noindent In particular, this proves Lemma \ref{lem:mono1}.\hfill$\Box$

\medskip\noindent{\bf Step \ref{sec:lim0}.4:} With the same kind of arguments as above, the following monotonicity holds (we omit the proof):
\begin{lem}\label{lem:mono2} 
Let $(\Vk)_{k\in\nn}\in C^0(B)$ and $(\uk)_{k\in\nn}\in C^4(B)$ such that (\ref{hyp:Vk}), (\ref{eq:uk}) and (\ref{bnd:nrj}) hold. We assume that $\uk$ is radially symmetrical for all $k\in\nn$. We assume that (\ref{max:uk}) holds. We assume that there exists $\delta_0\in (0,1)$ such that
$$\lim_{k\to +\infty}\mk^2\Delta\uk(\delta_0)=0.$$
We let $\delta\in (0,1)$ and $\eta\in (1,2)$. Then there exists $R_\eta>0$, there exists $(r_k)_{k\in\nn}\in \rr_{>0}$ such that $r_k\in [0,\delta]$ for all $k\in\nn$ and

\smallskip(i) $\lim_{k\to +\infty}\frac{\rk}{\mk}=+\infty$,\par
(ii) $r\mapsto r^\eta e^{\uk(r)}$ is decreasing on $[R_\eta\mk,\rk]$,\par
(iii) $\uk\to -\infty$ uniformly on $\overline{B}_\delta(0)\setminus B_{\rk}(0)$.
\end{lem}

\medskip\noindent{\bf Step \ref{sec:lim0}.4} We are in position to get the energy estimate for $e^{4\uk}$.

\begin{lem}\label{lem:mass0} 
Let $(\Vk)_{k\in\nn}\in C^0(B)$ and $(\uk)_{k\in\nn}\in C^4(B)$ such that (\ref{hyp:Vk}), (\ref{eq:uk}) and (\ref{bnd:nrj}) hold. We assume that $\uk$ is radially symmetrical for all $k\in\nn$. We assume that (\ref{max:uk}) holds. We assume that there exists $\delta_0\in (0,1)$ such that
$$\lim_{k\to +\infty}\mk^2\Delta\uk(\delta_0)=0.$$
Then for any $\delta\in (0,1)$, we have that:
$$\lim_{k\to +\infty}\int_{B_\delta(0)}\Vk  e^{4\uk}\, dx=16\pi^2.$$
In particular, $\Vk e^{4\uk}\rightharpoonup 16\pi^2\delta_0$ when $k\to +\infty$ in the sense of measures.
\end{lem}

\smallskip\noindent{\it Proof of Lemma \ref{lem:mass0}:} We prove the claim. We choose $\eta\in (1,2)$ and $R_\eta>0$, $(\rk)_{k\in\nn}$ as in Lemma \ref{lem:mono2}. We let $R>R_\eta$. It follows from Lemmae \ref{lem:cv0} and \ref{lem:mono2} that
\beq
&&\int_{B_\delta(0)\setminus B_{R\mk}(0)} e^{4\uk}\, dx \leq  \int_{B_{\rk}(0)\setminus B_{R\mk}(0)}\frac{(R\mk)^{4\eta}e^{4\uk(R\mk)}}{r^{4\eta}}\, dx+o(1)\\
&&\leq C(R\mk)^{4\eta}e^{4\uk(R\mk)}\int_{R\mk}^\delta r^{3-4\eta}\, dr+o(1)\leq C\frac{R^4\mk^4 e^{4\uk(R\mk)}}{\eta-1}+o(1)\\
&&\leq  \frac{C}{\eta-1}\left(\frac{\sqrt{96}R}{\sqrt{96}+R^2}\right)^4 +o(1)
\eeq
where $\lim_{k\to +\infty}o(1)=0$. Summing this integral and (\ref{lim:mass0}), letting $k\to +\infty$ and then $R\to +\infty$, we get the result. This proves Lemma \ref{lem:mass0}.
\hfill$\Box$

\medskip\noindent Point {\it (ii.a)} of Theorem \ref{thm:main} follows from Lemma \ref{lem:mass0}.


\begin{thebibliography}{1}

\bibitem{ars} Adimurthi; Robert, F.; Struwe, M. Concentration phenomena for Liouville equations in dimension four. {\it J. Eur. Math. Soc.}, to appear. Available on {\tt http://www-math.unice.fr/\~\,frobert}.

\bibitem{adistruwe} Adimurthi; Struwe, M. Global compactness properties of semilinear elliptic 
equations with critical exponential growth. {\it J. Funct. Analysis},
{\bf 175}, (2000), 125-167.

\bibitem{brezismerle} Br\'ezis, H.; Merle, F. Uniform estimates and blow-up behaviour for solutions of $-\Delta u=V(x)e^u$ in two dimensions. {\it Comm. Partial Differential Equations}, {\bf 16}, (1991),
1223-1253.

\bibitem{changchen}
Chang, S.-Y.A.; Chen, W. A note on a class of higher order conformally covariant equations. {\it Discrete and Continuous Dynamical systems}, {\bf 7}, (2001), 275-281.

\bibitem{druet} Druet, O. Multibumps analysis in dimension 2 - Quantification of blow up levels. {\it Duke Mathematical Journal}, to appear.

\bibitem{druetrobert} Druet, O.; Robert, F. Bubbling phenomena for fourth-order four-dimensional PDEs with exponential growth. {\it Proc. Amer. Math. Soc.}, {\bf 134}, (2006), 897-908. Available on {\tt http://www-math.unice.fr/\~\,frobert}.

\bibitem{gr} Ghoussoub, N.; Robert, F. Concentration estimates for Emden-Fowler equations with boundary singularities and critical growth. {\it IMRP}, to appear. Available on {\tt http://www-math.unice.fr/\~\,frobert}.

\bibitem{hebeyrobert} Hebey, E.; Robert, F. Coercivity and Struwe's compactness for Paneitz type operators with constant coefficients. {\it Calculus of Variations and Partial Differential Equations}, {\bf 13}, (2001), 491-517.

\bibitem{hrw} Hebey, E.; Robert, F.; Wen, Y. Compactness and global estimates for a fourth 
order equation of critical Sobolev growth arising from conformal geometry, {\it Communication in Contemporary Mathematics}, to appear.

\bibitem{lishafrir} Li, Y.; Shafrir, I. Blow-up analysis for solutions of $-\Delta u=Ve^u$ in dimension two. {\it Indiana Univ. Math. J.}, {\bf 43}, (1994), 1255-1270.

\bibitem{lin} Lin, C.S. A classification of solutions of a conformally invariant fourth 
order equation in ${\mathbb R}^n$. {\it Comment. Math. Helv.}, {\bf 73}, 1998, 206-231.

\bibitem{malchiodi} Malchiodi, A. Compactness of solutions to some
  geometric fourth-order equations. {\it J. Reine Angew. Math.}, to appear.

\bibitem{malchiodistruwe} Malchiodi, A.; Struwe, M. The $Q$-curvature
  flow on $S^4$. {\it Preprint 2004}.


\bibitem{r1} Robert, F. Quantization effects for a fourth order equation of exponential growth in dimension four. {\it Preprint 2005}. Available on {\tt http://www-math.unice.fr/\~\,frobert}.

\bibitem{robertstruwe}
Robert, F.; Struwe, M.  Asymptotic profile for a fourth order pde with critical exponential growth in dimension four. {\it Advanced Nonlinear Studies}, {\bf 4}, (2004), 397-415.

\bibitem{schoenzhang}
Schoen, R.; Zhang, D. Prescribed scalar curvature on the $n$-sphere. {\it Calc. Var. Partial Differential Equations}, {\bf 4}, (1996), 1-25. 

\bibitem{t} Tarantello, G. A quantization property for blow-up solutions of singular Liouville-type equations. {\it J. Funct. Anal.}, {\bf 219}, (2005), 368-399.

\bibitem{wei} Wei, J. Asymptotic behavior of a nonlinear fourth order eigenvalue problem. {\it  Comm. Partial Differential Equations}, {\bf 21}, (1996),  no. 9-10, 1451-1467.


\end{thebibliography}
\end{document}